# Attractiveness of the Haar measure for linear cellular automata on Markov subgroups

Alejandro Maass[1,*], Servet Martínez[1,*],
Marcus Pivato[2,†] and Reem Yassawi[2,†]

*Universidad de Chile and Trent University*

**Abstract:** For the action of an algebraic cellular automaton on a Markov subgroup, we show that the Cesàro mean of the iterates of a Markov measure converges to the Haar measure. This is proven by using the combinatorics of the binomial coefficients on the regenerative construction of the Markov measure.

## 1. Introduction and main concepts

Let $(\mathcal{A}, +)$ be a finite Abelian group with unit 0 and $(\mathcal{A}^{\mathbb{Z}}, +)$ be the product group with the componentwise addition. The shift map $\sigma : \mathcal{A}^{\mathbb{Z}} \to \mathcal{A}^{\mathbb{Z}}$ is defined on $\mathbf{x} \in \mathcal{A}^{\mathbb{Z}}$ and $n \in \mathbb{Z}$ by $(\sigma\mathbf{x})_n = x_{n+1}$.

Let $\mu$ be a probability measure on $\mathcal{A}^{\mathbb{Z}}$ and $\Phi = \sum_{i=\ell}^{r} c_i \sigma^i$ a linear cellular automaton on $\mathcal{A}^{\mathbb{Z}}$, where $\ell, r, c_\ell, \ldots, c_r$ are integers such that $\ell \leq r$. The study of the convergence of the sequence $(\Phi^n(\mu) : n \in \mathbb{N})$ started with the pioneering work of Lind in [8]. This work, devoted to the linear cellular automaton $\Phi = \sigma^{-1} + id + \sigma$ on $\{0,1\}^{\mathbb{Z}}$ and $\mu$ a Bernoulli measure, stated that the Cesàro mean of the sequence converges to the uniform Bernoulli measure. So $\Phi$ randomizes any Bernoulli measure.

After such work two different strategies were developed to prove the same kind of results but for more general classes of linear cellular automata and other kind of starting measures. The combinatorial properties of the Pascal triangle and the technique of regeneration of measures proposed in [1] were exploited in [2, 3, 5, 9, 10] to prove Lind's result in the case $\Phi = id + \sigma$ and $\mu$ is a measure with complete connections and summable memory decay for any finite group $\mathcal{A}$. In [12, 13] an harmonic analysis formalism was developed to prove that *diffusive* linear cellular automata randomize *harmonically mixing* measures (both concepts were defined in [12]). Other related results were proved in [5] and [14]

In this paper we consider the same problem for linear cellular automata defined in a Markov subgroup. Before state our main result we need some precise notations and definitions. More details and references on the subject can be found in [6, 7, 15].

---

*Supported by Nucleus Millennium Information and Randomness P04-069-F.
†Partially supported by the NSERC Canada.
[1]Departamento de Ingeniería Matemática and Centro de Modelamiento Matemático, Universidad de Chile, Av. Blanco Encalada 2120 5to piso, Santiago, Chile, e-mail: `amaass@dim.uchile.cl` e-mail: `smartine@dim.uchile.cl`
[2]Department of Mathematics, Trent University, Canada, e-mail: `pivato@xaravve.trentu.ca` e-mail: `ryassawi@xaravve.trentu.ca`

Let $\mathfrak{G} \subseteq \mathcal{A}^{\mathbb{Z}}$ be a Markov subgroup, that is, it is a topological Markov shift defined by a $0-1$ incidence matrix $\mathbf{M}$ and it is a subgroup of $\mathcal{A}^{\mathbb{Z}}$. For every $g \in \mathcal{A}$, let $\mathcal{F}(g) = \{h \in \mathcal{A} : \mathbf{M}_{gh} = 1\}$ be the set of followers of $g$ and $\mathcal{P}(g) = \{h \in \mathcal{A} : \mathbf{M}_{hg} = 1\}$ be the set of predecessors of $g$. Then $\mathcal{F} := \mathcal{F}(0)$ and $\mathcal{P} := \mathcal{P}(0)$ are subgroups of $\mathcal{A}$ and $\mathcal{F}(g)$ and $\mathcal{P}(g)$ are cosets of $\mathcal{F}$ and $\mathcal{P}$ in $\mathcal{A}$ respectively. It holds that $|\mathcal{F}| = |\mathcal{P}|$, where $|\cdot|$ means the cardinality of a set. We fix a function $f : \mathcal{A} \to \mathcal{A}$ such that $f(g) \in \mathcal{F}(g)$ for every $g \in \mathcal{A}$, and then $\mathcal{F}(g) = f(g) + \mathcal{F}$.

For every $g \in \mathcal{A}$ and $n \geq 1$ define $\mathcal{F}^n(g) = \{h \in \mathcal{A} : \mathbf{M}^n_{gh} > 0\}$. One has $\mathcal{F}^{n+1}(g) = \cup_{h \in \mathcal{F}(g)} \mathcal{F}^n(h)$ and $|\mathcal{F}^n(g)| = |\mathcal{F}^n(0)|$ for every $g \in \mathcal{A}$. The Markov shift $\mathfrak{G}$ is transitive if and only if $\mathbf{M}$ is irreducible, that is, there exists $n \geq 1$ such that $\mathcal{F}^n(0) = \mathcal{A}$. Denote by $r$ the smallest $n$ verifying this condition. In this case we have the mixing property $\mathcal{F}^n(g) = \mathcal{A}$ for every $n \geq r$ and $g \in \mathcal{A}$.

For $n \geq 0$ and $g_0, g_n \in \mathcal{A}$ such that $\mathbf{M}^n_{g_0 g_n} > 0$ define

$$\mathfrak{C}^n(g_0, g_n) = \{(g_1, \ldots, g_{n-1}) \in \mathcal{A}^{n-1} : \mathbf{M}_{g_i g_{i+1}} = 1, i \in \{0, \ldots, n-1\}\}.$$

Then $|\mathfrak{C}^n(g_0, g_n)| = |\mathfrak{C}^n(0, 0)|$. Therefore, in the transitive case, $|\mathfrak{C}^n(g, h)| = |\mathfrak{C}^n(0, 0)|$ for every $g, h \in \mathcal{A}$ and $n \geq r$.

The Haar measure of $(\mathfrak{G}, +)$ is denoted by $\nu$. It is the Markov measure defined by the stochastic matrix $\mathbf{L} = (L_{gh} : g, h \in \mathcal{A})$ with $L_{gh} = |\mathcal{F}|^{-1}$ if $h \in \mathcal{F}(g)$ and $L_{gh} = 0$ otherwise, and the $\mathbf{L}$–stationary vector $\rho = (\rho(g) = |\mathcal{A}|^{-1} : g \in \mathcal{A})$. The Haar measure is the maximal entropy measure for the Markov shift $(\mathfrak{G}, \sigma)$.

It is useful to introduce a notation for the finite-dimensional distributions associated to $\nu$. Thus, for $\ell \geq 1$, $\nu^{(\ell)}$ is the finite-dimensional distribution of $\nu$ in $\mathcal{A}^\ell$, and it is concentrated in the subset $\mathfrak{G}_\ell = \{(g_0, \ldots, g_{\ell-1}) \in \mathcal{A}^\ell : g_{i+1} \in \mathcal{F}(g_i), \forall i \in \{0, \ldots, \ell-2\}\}$. Hence,

$$\text{if } \gamma_\ell = \frac{1}{|\mathcal{A}| \cdot |\mathcal{F}|^{\ell-1}}, \quad \text{then} \quad \nu^{(\ell)}(\mathbf{g}) = \begin{cases} \gamma_\ell & \text{if } \mathbf{g} \in \mathfrak{G}_\ell, \\ 0 & \text{if } \mathbf{g} \notin \mathfrak{G}_\ell. \end{cases}$$

### 1.1. Main Result

Let $\mu$ be a Markov probability measure in $\mathcal{A}^{\mathbb{Z}}$ that is defined by a stochastic matrix $\mathbf{P}$ and by some probability vector $\pi$ invariant for $\mathbf{P}$ ($\pi \cdot \mathbf{P} = \pi$). The measure $\mu$ is said to be compatible with $\mathbf{M}$ if $\mathbf{P}_{gh} > 0$ if and only if $\mathbf{M}_{gh} > 0$. This property is equivalent to the fact that the support of $\mu$ is $\mathfrak{G}$.

If $\mathbf{x} \in \mathcal{A}^{\mathbb{Z}}$ and $m, n \in \mathbb{Z}$ are two integers with $m \leq n$, define

$$\mathbf{x}_m^n = (x_m, x_{m+1}, \ldots, x_n).$$

Consider the endomorphism $\Phi = id + \sigma$ on $\mathcal{A}^{\mathbb{Z}}$, that is $(\Phi \mathbf{x})_n = x_n + x_{n+1}$ for $\mathbf{x} \in \mathcal{A}^{\mathbb{Z}}$ and $n \in \mathbb{Z}$. Our main result is the following one.

**Theorem 1.1.** *Assume $\mu$ is a Markov probability measure in $\mathcal{A}^{\mathbb{Z}}$ compatible with $\mathbf{M}$. Furthermore, assume the Abelian group $\mathcal{A}$ is $p^s$-torsion for some prime number $p$ and $s \geq 1$ and that $\mathbf{M}$ is irreducible. Then the Cesàro mean of $\mu$, under the action of $\Phi$, converges to the Haar measure $\nu$. That is, for any $m \in \mathbb{N}$ and $\mathbf{g} \in \mathfrak{G}_m$:*

$$\lim_{N \to \infty} \frac{1}{N} \sum_{n=0}^{N-1} \mu\left( (\Phi^n \mathbf{x})_0^{m-1} = \mathbf{g} \right) = \nu(\mathbf{g}). \tag{1.1}$$

The elements of our proof includes a regenerative construction of the Markov measure and the combinatorics of the binomial coefficients associated with the iterates of the cellular automaton. In [11] a more general result is shown: for Markov fields verifying a 'filling property' it is shown that the attractive property of the Theorem holds. In the one-dimensional case this 'filling property' is always satisfied, therefore the result follows.

## 2. Construction of a Markov measure and renewal properties

We use the Athreya–Ney [1] representation of Markov chains. It says that it is possible to enlarge the probability space we are considering in order to include a family of integer random times.

### 2.1. Construction of a Markov measure

Let $\mathfrak{G} \subseteq \mathcal{A}^{\mathbb{Z}}$ be a Markov subgroup with incidence matrix $\mathbf{M}$. Let $\mu$ be a Markov probability measure in $\mathcal{A}^{\mathbb{Z}}$ compatible with $\mathbf{M}$ defined by a pair $(\pi, \mathbf{P})$ as in subsection 1.1. We describe a procedure to construct the restriction of $\mu$ to $\mathcal{A}^{\mathbb{N}}$. In this purpose it is useful to introduce the following notation: for $g \in \mathcal{A}$ we put $\mu_g$ the measure induced on $\mathcal{A}^{\mathbb{N}}$ by $\mu$ conditioned to the event $\{x_{-1} = g\}$; then $\sum_{g \in \mathcal{A}} \pi(g)\mu_g$ coincides with the restriction of $\mu$ to $\mathcal{A}^{\mathbb{N}}$.

Let $\alpha > 0$ be a strictly positive number such that

$$\alpha < \min\{\mathbf{P}_{gh} : g, h \in \mathcal{A}, \mathbf{M}_{gh} > 0\}.$$

We consider a probability space $(\Omega, \mathcal{B}, \mathbb{P})$ and three independent processes of i.i.d. random variables defined on this space $\mathbf{U} = (U_n : n \in \mathbb{N})$, $\mathbf{W} = (W_n : n \in \mathbb{N})$ and $\mathbf{V} = (V_n : n \in \mathbb{N})$ whose marginal distributions are as follows: $U_n$ is Bernoulli($\alpha$), that is $\mathbb{P}(U_n = 1) = \alpha = 1 - \mathbb{P}(U_n = 0)$; $W_n$ is uniformly distributed in $\mathcal{F}$, so $\mathbb{P}(W_n = g) = |\mathcal{F}|^{-1}$ for $g \in \mathcal{F}$; and $V_n$ is uniformly distributed in the unit interval $[0, 1]$.

Let us construct a sequence $(x_n : n \in \mathbb{N}) \in \mathcal{A}^{\mathbb{N}}$ as a deterministic function of $(\mathbf{U}, \mathbf{W}, \mathbf{V})$. For any $g, h \in \mathcal{A}$ define

$$\mathbf{Q}_{gh} = \frac{\mathbf{P}_{gh} - \alpha |\mathcal{F}|^{-1} \mathbf{M}_{gh}}{1 - \alpha}.$$

Thus, $\mathbf{Q} = (\mathbf{Q}_{gh} : g, h \in \mathcal{A})$ is a stochastic matrix compatible with $\mathbf{M}$. For each $g \in \mathcal{A}$ fix $\{\tilde{Q}_{gh} \subseteq [0, 1] : h \in \mathcal{F}(g)\}$ a measurable partition of the interval $[0, 1]$ such that the Lebesgue measure of $\tilde{Q}_{gh}$ is $\mathbf{Q}_{gh}$. For all $g \in \mathcal{A}$, $u \in \{0, 1\}$, $w \in \mathcal{F}$ and $v \in [0, 1]$ define

$$H(g, u, w, v) = u\left(f(g) + w\right) + (1-u) \sum_{h \in \mathcal{F}(g)} h \, \mathbb{1}_{\widetilde{Q}_{gh}}(v).$$

Now, for any $x_{-1} \in \mathcal{A}$ and for $n \geq 0$ we put

$$x_n = H(x_{n-1}, U_n, W_n, V_n).$$

It is clear that the distribution of the sequence $\mathbf{x} = (x_n : n \in \mathbb{N})$ is $\mu_{x_{-1}}$. If $x_{-1} \in \mathcal{A}$ is a random variable with distribution $\pi$, then the distribution of $\mathbf{x}$ is $\mu$.

## 2.2. The associated renewal process

For any $s, t \in \mathbb{N}$ with $s \leq t$ define

$$\left( \mathbf{U}_s^t = \mathbf{1} \right) \iff \left( U_k = 1, \text{ for all } k \in \{s, ..., t\} \right).$$

For every $m \geq 1$ define a renewal process $(T_n^{(m)} : n \in \mathbb{N})$ given by:

$$T_0^{(m)} = 0, \ T_1^{(m)} = \min\{i > T_0^{(m)} : \mathbf{U}_i^{i+m} = \mathbf{1}\},$$

and

$$T_n^{(m)} = \min\{i > T_{n-1}^{(m)} + m : \mathbf{U}_i^{i+m} = \mathbf{1}\} \text{ for } n \geq 2.$$

It is clear that $(T_{n+1}^{(m)} - T_n^{(m)} : n \geq 1)$ is a family of i.i.d. random variables. Also, from our computations below, it follows that the distribution of $T_1^{(m)}$ has a geometrical tail; that is, there exists $\beta := \beta(m) \in [0, 1)$ such that $\mathbb{P}(T_1^{(m)} > t) \leq \beta^t$ for any $t \geq 0$.

Let $\mathbf{N}^{(m)}$ be the renewal process induced by $(T_n^{(m)} : n \in \mathbb{N})$. That is, for any $A \subseteq \mathbb{N}$,

$$\mathbf{N}^{(m)}(A) = \{n \in A : \text{ for some } \ell \in \mathbb{N}, \ T_\ell^{(m)} = n\}.$$

Let $n \in A$ with $n > 0$. One has that $U_{n-1} = 0$ and $\mathbf{U}_n^{n+m} = \mathbf{1}$ implies $n \in \mathbf{N}^{(m)}(A)$. Then $\mathbb{P}(n \in \mathbf{N}^{(m)}(A)) \geq \alpha^{m+1}(1 - \alpha) := \delta > 0$. Clearly if $0 \in A$ then $\mathbb{P}(0 \in \mathbf{N}^{(m)}(A)) = 1 \geq \delta$.

For $A$ a finite subset of $\mathbb{N}$ and $m \in \mathbb{N}$ we say $A$ is *m-separated* if, for any $a, b \in A$ with $a \neq b$, $|a - b| \geq m + 1$. One gets,

$$\left( A \text{ is } m\text{-separated} \right) \implies \left( \mathbb{P}\left(\mathbf{N}^{(m)}(A) \neq \emptyset\right) \geq 1 - (1-\delta)^{|A|} \right). \tag{2.1}$$

Let $A^{(m)}$ be the largest $m$-separated subset of $A$; then $|A^{(m)}| \geq |A|/(m+1)$. Thus,

$$\mathbb{P}\left( \mathbf{N}^{(m)}(A) = \emptyset \right) \leq \mathbb{P}\left( \mathbf{N}^{(m)}(A^{(m)}) = \emptyset \right) \leq \left((1-\delta)^{1/(m+1)}\right)^{|A|}.$$

Hence, the distribution of $T_1^{(m)}$ has a geometric tail.

## 3. Convergence of the Cesàro mean

### 3.1. Independence lemmas

Assume $\mathfrak{G}$ is transitive and recall $r$ is the smallest integer verifying $\mathcal{F}^r(0) = \mathcal{A}$. Also recall that $\gamma_\ell = |\mathcal{A}|^{-1}|\mathcal{F}|^{-(\ell-1)}$. Let $m \leq n$ in $\mathbb{Z}$. If $\mathbf{x}$ is a random variable in $\mathcal{A}^\mathbb{Z}$ with distribution $\mu$ define $\mathcal{F}_m^n$ to be the sigma-algebra generated by $\mathbf{x}_m^n$.

**Lemma 3.1.** *Let $k \geq r$ and $m \geq 0$. Then the random variable $\mathbf{x}_k^{k+m}$ conditioned to $\left(\mathbf{U}_{k-r}^{k+m} = \mathbf{1}, \mathcal{F}_0^{k-r-1}\right)$ is $\nu^{(m+1)}$–distributed. That is, for any $\mathbf{g} \in \mathfrak{G}_{m+1}$,*

$$\mathbb{P}\left( \mathbf{x}_k^{k+m} = \mathbf{g} \ \Big| \ \mathbf{U}_{k-r}^{k+m} = \mathbf{1}, \ \mathcal{F}_0^{k-r-1} \right) = \gamma_{m+1}.$$

*Also, the variable $\mathbf{x}_k^{k+m}$ conditioned to $(\mathbf{U}_{k-r}^{k+m+r} = \mathbf{1}, \mathcal{F}_0^{k-r-1} \vee \mathcal{F}_{k+m+r+1}^n)$ is $\nu^{(m+1)}$–distributed for any $n \geq k + m + r + 1$. That is, for any $\mathbf{g} \in \mathfrak{G}_{m+1}$,*

$$\mathbb{P}\left(\mathbf{x}_k^{k+m} = \mathbf{g} \ \Big| \ \mathbf{U}_{k-r}^{k+m+r} = \mathbf{1}, \ \mathcal{F}_0^{k-r-1} \vee \mathcal{F}_{k+m+r+1}^n\right) = \gamma_{m+1}.$$

*Proof.* Let $\mathbf{g} = (g_0, \ldots, g_m) \in \mathfrak{G}_{m+1}$ and put $n = k - r$. For any fixed $\mathbf{h} = (h_0, \ldots, h_{n-1}) \in \mathfrak{G}_n$,

$$\mathbb{P}\left(\mathbf{x}_k^{k+m} = \mathbf{g} \middle| \mathbf{x}^{n-1} = \mathbf{h} \text{ and } \mathbf{U}_n^{k+m} = \mathbf{1}\right)$$

$$= \sum_{\mathbf{z} \in \mathfrak{C}^{r+1}(h_{n-1}, g_0)} \mathbb{P}\left(\mathbf{x}_n^{n-1} = \mathbf{z} \text{ and } \mathbf{x}^{k+m} = \mathbf{g} \middle| \mathbf{x}_0^{n-1} = \mathbf{h} \text{ and } \mathbf{U}_n^{k+m} = \mathbf{1}\right)$$

$$= \sum_{\mathbf{z} \in \mathfrak{C}^{r+1}(h_{n-1}, g_0)} |\mathcal{F}|^{-(r+m+1)} = |\mathfrak{C}^{r+1}(0,0)| \cdot |\mathcal{F}|^{(r+m+1)}$$

$$= \frac{|\mathcal{F}|^{r+1}}{|\mathcal{A}|} \cdot |\mathcal{F}|^{-(r+m+1)} = \frac{1}{|\mathcal{A}| \cdot |\mathcal{F}^m|} = \gamma_{m+1}$$

This proves the first part, the second one is entirely analogous. $\square$

Now assume that $\mathcal{A}$ is $p^s$-torsion for some prime number $p$ and some $s \geq 1$, with $s$ being the smallest number verifying this property. That is: $mg = 0$ for every $g \in \mathcal{A}$ and $m \in p^s\mathbb{Z}$, and for every $m < p^s$ there exists some $g \in \mathcal{A}$ such that $mg \neq g$. Observe that for every $c \in \mathbb{Z}_{p^s}$ relatively prime to $p$, there exists a multiplicative inverse $c^{-1} \in \mathbb{Z}_{p^s}$, such that $cc^{-1} = 1 \mod (p^s)$ and $c^{-1}$ is also relatively prime to $p$. Thus, for any $g \in \mathcal{A}$, $cc^{-1}g = g$. Moreover, $(cg = h) \iff (g = c^{-1}h)$.

### 3.2. The transformation

Recall $\sigma$ is the shift map in $\mathcal{A}^{\mathbb{Z}}$ and $\Phi = id + \sigma$ is an endomorphism of $\mathcal{A}^{\mathbb{Z}}$. Fix $\mathbf{x} \in \mathcal{A}^{\mathbb{Z}}$. Then for all $n \geq 0$ and $i \in \mathbb{Z}$, one has

$$(\Phi^n \mathbf{x})_i = \sum_{k=0}^n \binom{n}{k} x_{i+k}. \tag{3.1}$$

For every $m \in \mathbb{N}$ denote by $m^{(s)}$ its equivalent class $\mod (p^s)$ in $\mathbb{Z}_{p^s}$. Hence, for all $n \geq 0$ and $i \in \mathbb{Z}$, equality (3.1) can be written

$$(\Phi^n \mathbf{x})_i = \sum_{k=0}^n \binom{n}{k}^{(s)} x_{i+k}. \tag{3.2}$$

Let $m, \ell \geq 0$. For $n \geq 0$ and $0 \leq k \leq n$ we say $k$ is $(m, \ell)$–*isolated* in $n$ if $\binom{n}{k}^{(s)}$ is relatively prime to $p$, while, if $k' \in \{k-m, \ldots, k+\ell\}$ with $k' \neq k$ then $\binom{n}{k'}^{(s)} = 0$. Here the convention is $\binom{n}{k'}^{(s)} = 0$ whenever $k' < 0$ or $k' > n$.

**Lemma 3.2.** *Let $m \in \mathbb{N}$ and $n \geq 2r+2m+1$. If $0 \leq k \leq n$ is $(r+m, r+m)$-isolated in $n$, then for every $i \in \mathbb{Z}$ and $\mathbf{g} \in \mathfrak{G}_{m+1}$, one has*

$$\mathbb{P}\left((\Phi^n \mathbf{x})_i^{i+m} = \mathbf{g} \,\middle|\, \mathbf{U}_{i+k-r}^{i+k+r+m} = \mathbf{1}\right) = \gamma_{m+1}.$$

*Proof.* Since $\mu$ is $\sigma$-invariant then $\Phi^n(\mu)$ is $\sigma$-invariant too; hence it suffices to prove the result for $i = 0$. In other words, it suffices to show that

$$\mathbb{P}\left((\Phi^n \mathbf{x})_0^m = \mathbf{g} \,\middle|\, \mathbf{U}_{k-r}^{k+r+m} = \mathbf{1}\right) = \gamma_{m+1}.$$

Consider $j \in \{0, ..., m\}$. Define $Y_j = \binom{n}{k}^{(s)} x_{j+k}$ and

$$X_j = \sum_{\substack{k'=0 \\ k' \neq k}}^{n} \binom{n}{k'}^{(s)} x_{j+k'} \overset{(*)}{=} \sum_{k'=0}^{k-r-m-1} \binom{n}{k'}^{(s)} x_{j+k'} + \sum_{k'=k+r+m+1}^{n} \binom{n}{k'}^{(s)} x_{j+k'}, \tag{3.3}$$

where $(*)$ is because $k$ is $(r+m, r+m)$-isolated.

Thus, $(\Phi^n \mathbf{x})_j = X_j + Y_j$. If $\mathbf{X} = (X_0, \ldots, X_m)$ and $\mathbf{Y} = (Y_0, \ldots, Y_m)$, then $(\Phi^n \mathbf{x})_0^m = \mathbf{X} + \mathbf{Y}$. One gets,

$$\mathbb{P}\left((\Phi^n \mathbf{x})_0^m = \mathbf{g} \middle| \mathbf{U}_{k-r}^{k+r+m} = \mathbf{1}\right) = \mathbb{P}\left(\mathbf{X} + \mathbf{Y} = \mathbf{g} \middle| \mathbf{U}_{k-r}^{k+r+m} = \mathbf{1}\right)$$
$$= \sum_{\mathbf{h} \in \mathcal{A}^{m+1}} \mathbb{P}\left(\mathbf{Y} = \mathbf{g} - \mathbf{h} \middle| \mathbf{X} = \mathbf{h}, \mathbf{U}_{k-r}^{k+r+m} = \mathbf{1}\right)$$
$$\times \mathbb{P}\left(\mathbf{X} = \mathbf{h} \middle| \mathbf{U}_{k-r}^{k+r+m} = \mathbf{1}\right). \tag{3.4}$$

Let $c = \binom{n}{k}^{(s)}$; then $c$ is relatively prime to $p$, and $\mathbf{Y} = c \cdot \mathbf{x}_k^{k+m}$. Thus, if $c^{-1}$ is the (mod $p^s$) inverse of $c$, then $\mathbf{x}_k^{k+m} = c^{-1} \cdot \mathbf{Y}$. Thus, for any fixed $\mathbf{h} \in \mathcal{A}^{m+1}$,

$$\mathbb{P}\left(\mathbf{Y} = \mathbf{g} - \mathbf{h} \middle| \mathbf{X} = \mathbf{h} \text{ and } \mathbf{U}_{k-r}^{k+r+m} = \mathbf{1}\right) \tag{3.5}$$
$$= \mathbb{P}\left(\mathbf{x}_k^{k+m} = c^{-1}(\mathbf{g} - \mathbf{h}) \middle| \mathbf{X} = \mathbf{h} \text{ and } \mathbf{U}_{k-r}^{k+r+m} = \mathbf{1}\right)$$
$$\overset{(\dagger)}{=} \gamma_{m+1},$$

where $(\dagger)$ is because equation (3.3) implies that $\mathbf{X}$ is a function only of $\mathbf{x}_0^{k-r-1}$ and $\mathbf{x}_{k+r+m+1}^{n+m}$; and it allows to apply Lemma 3.1. Substituting (3.5) into (3.4) yields

$$\mathbb{P}\left((\Phi^n \mathbf{x})_0^m = \mathbf{g} \middle| \mathbf{U}_{k-r}^{k+r+m} = \mathbf{1}\right) = \sum_{\mathbf{h} \in \mathcal{A}^{m+1}} \gamma_{m+1} \cdot \mathbb{P}\left(\mathbf{X} = \mathbf{h} \middle| \mathbf{U}_{k-r}^{k+r+m} = \mathbf{1}\right)$$
$$= \gamma_{m+1}.$$

Therefore the result follows. □

### 3.3. Elementary facts on the Pascal triangle

Now we use the following result on the Pascal triangle. Let $n = \sum_{j \in \mathbb{N}} n_j p^j$ be the decomposition of $n$ in base $p$, so $n_j \in \{0, \ldots, p-1\}$ for every $j \in \mathbb{N}$. For $0 \leq k \leq n$ consider the decompositions $k = \sum_{j \in \mathbb{N}} k_j p^j$ and $n - k = \sum_{j \in \mathbb{N}} (n-k)_j p^j$. The Kummer's Theorem on binomial coefficients, whose proof can be found in [4], states the following result.

**Lemma 3.3.** *The biggest integer $\ell$ such that $p^\ell$ divides $\binom{n}{k}$ is the number of carries needed to sum $k$ and $n-k$ in base $p$.*

We introduce the following notation. For $n \in \mathbb{N}$ and $i \geq 0$ we put

$$J_i(n) = \{j \geq i \ : \ n_j \neq 0\}, \quad \xi_i(n) = |J_i(n)|.$$

For a real number $c$ denote by $\lfloor c \rfloor$ the integer part of $c$ and for $a \geq 0$ define

$$p^a \mathbb{N} = \{p^a \cdot n \ : \ n \in \mathbb{N}\} = \{n \in \mathbb{N} \ : \ n_j = 0 \text{ for all } 0 \leq j < a\}.$$

**Lemma 3.4.** *Let $m \geq 1$ and $a \geq 2s+1$ be such that $p^{\lfloor a/2 \rfloor} > m$. For every $n \in p^a \mathbb{N}$ and $i \geq a$,* $\left|\{0 \leq k \leq n : k \text{ is } (m,m)\text{-isolated in } n\}\right| \geq 2^{\xi_{a+i}(n)} - 1.$

*Proof.* Fix a nonempty subset $J \subseteq J_{i+a}(n)$, and define $0 \leq k \leq n$ by $k_j = n_j$ if $j \notin J$ and $k_j = n_j - 1$ for $j \in J$. Therefore $n - k$ verifies $(n-k)_j = 1$ for $j \in J$ and $(n-k)_j = 0$ for $j \notin J$. Then there is no carry in the sum of $k$ and $n-k$, so Lemma 3.3 says that $\binom{n}{k}$ is relatively prime to $p$. It remains to show that $\binom{n}{k'}^{(s)} = 0$ for all $k' \in \{k-m, \ldots, k+m\} \setminus \{k\}$.

**Case 1:** Let $1 \leq v \leq m$ and $k' = k - v$. Then the $p$-ary decomposition of $k'$ has some nonzero elements in coordinates between $0$ and $a-1$ (because $p^a$ divides $k$, but does not divide $v$); moreover, it has at least $b = \lfloor a/2 \rfloor$ zeros in $\{0, \ldots, a-1\}$. However, the $p$-ary decomposition of $n$ has no nonzero elements in $\{0, \ldots, a-1\}$, so there must be at least $b$ carries in the addition: $(n-k') + k' = n$. Thus, Lemma 3.3 says that $p^b$ divides $\binom{n}{k'}$. But $b = \lfloor a/2 \rfloor \geq \lfloor (2s+1)/2 \rfloor \geq s$, so we conclude that $\binom{n}{k'}^{(s)} = 0$.

**Case 2:** Let $k' = k + v$ with $1 \leq v \leq m$. Then $k'_j \geq n_j$ for every $j < a+i$ and for some $j' < a$ we have $k'_{j'} > n_{j'}$. Hence, the sum in base $p$ of $k'$ and $n-k'$ will have at least $a$ carries, so Lemma 3.3 says that $p^a$ divides $\binom{n}{k'}$ and finally $\binom{n}{k'}^{(s)} = 0$. □

**Lemma 3.5.** *Let $a \geq 0$. Then, the set*

$$\mathcal{M}_0 = \left\{ n \in p^a \mathbb{N} \ : \ \xi_{a + \lfloor \frac{1}{2} \log_p(n) \rfloor}(n) \geq \frac{1}{5} \log_p(n) \right\}$$

*is of density $1$ in $p^a \mathbb{N}$. That is,*

$$\lim_{N \to \infty} \frac{|\mathcal{M}_0 \cap \{0, \ldots, N-1\}|}{|p^a \mathbb{N} \cap \{0, \ldots, N-1\}|} = 1.$$

*Proof.* Let $\mathcal{M} = \left\{ n \in \mathbb{N} \ : \ \xi_{\lfloor \frac{1}{2} \log_p(n) \rfloor}(n) \geq \frac{1}{5}(\log_p(n) + a) \right\}$. Then

$$\lim_{N \to \infty} \frac{1}{N} \left| \mathcal{M} \cap \{0, \ldots, N-1\} \right| = 1.$$

To see this, let $n \in \mathbb{N}$ be a 'generic' large integer; then the Law of Large Numbers says that only about $\frac{1}{p}$ of the $p$-ary digits of $n$ are zero; hence $\frac{p-1}{p}$ are nonzero. Since there are $\frac{1}{2} \log_p(n)$ digits in the range $[\frac{1}{2} \log_p(n), \log_p(n)]$, we conclude, with asymptotic probability $1$, that at least $\frac{p-1}{2p} \log_p(n)$ digits in $[\frac{1}{2} \log_p(n), \log_p(n)]$ are nonzero; hence $\xi_{\lfloor \frac{1}{2} \log_p(n) \rfloor}(n) \geq \frac{p-1}{2p} \log_p(n) \geq \frac{1}{5}(\log_p(n) + a)$ (assuming $p \geq 2$ and $\log_p(n) > 4a$).

Now define bijection $\psi : \mathbb{N} \to p^a \mathbb{N}$ by $\psi(n) = p^a n$; then $\psi(\mathcal{M}) \supseteq \mathcal{M}_0$. The lemma follows. □

### 3.4. Proof of Theorem 1.1

*Proof.* Fix $m \geq 1$. If $\mathcal{N} \subset \mathbb{N}$, we say that the $m$-dimensional marginal of the Cesàro mean converges along $\mathcal{N}$ if for any $\mathbf{g} \in \mathfrak{G}_m$

$$\lim_{N \to \infty} \frac{1}{|\mathcal{N} \cap \{0, \ldots, N-1\}|} \sum_{n \in \mathcal{N} \cap \{0, \ldots, N-1\}} \mu\left( (\Phi^n \mathbf{x})_0^{m-1} = \mathbf{g} \right) = \nu(\mathbf{g}). \quad (3.6)$$

Let $m' = m + r$ and consider $a \geq 2s + 1$ with $p^{\lfloor a/2 \rfloor} > m'$ (as in Lemma 3.4). Let $\mathcal{M}_0$ be as in Lemma 3.5. We claim that the $m$-dimensional marginal of the Cesàro mean converges along $\mathcal{M}_0$.

Let $n \in \mathcal{M}_0$ be enough large such that $i = \lfloor \frac{1}{2} \log_p(n) \rfloor \geq a$. Define $A = \{0 \leq k \leq n : k \text{ is } (m', m')\text{-isolated in } n \}$. Therefore, by Lemma 3.4 and the definition of $\mathcal{M}_0$,

$$|A| \geq 2^{\xi_{a+i}(n)} - 1 \geq 2^{\frac{1}{5} \log_p(n)} - 1 = 2^{\frac{1}{C} \log_2(n)} - 1 = n^{1/C} - 1, \qquad (3.7)$$

where $C = 5 \log_2(p)$. Thus,

$$\mu\left(\exists k \in A \text{ with } \mathbf{U}_{k-r}^{k+m'} = \mathbf{1}\right) \geq \mu\left(\mathbf{N}^{(2m')}(A - m') \neq \emptyset\right)$$
$$\underset{(2.1)}{\geq} 1(1 - \delta)^{|A|} \underset{(3.7)}{\geq} 1 - (1 - \delta)^{n^{1/c} - 1},$$

where "(3.7)" is by equation (3.7) and "(2.1)" is by equation (2.1) (since $A$ is $(2m')$-separated).

Finally, Lemma 3.2 implies $|\mu((\Phi^n \mathbf{x})_i^{i+m-1} = \mathbf{g}) - \gamma_m| \leq (1-\delta)^{n^{1/C}-1}$, and thus, $\lim_{n \to \infty, n \in \mathcal{M}_0} |\mu((\Phi^n \mathbf{x})_i^{i+m-1} = \mathbf{g}) - \gamma_m| = \lim_{n \to \infty} (1-\delta)^{n^{1/C}-1} = 0$, as desired.

Lemma 3.5 then implies that the $m$-dimensional marginal of the Cesàro mean converges along $p^a \mathbb{N}$. Now, since $\nu$ is invariant for powers of $\Phi$, we find that for any $0 \leq j < p^a$, the $m$-dimensional marginal of the Cesàro mean also converges along $\mathcal{M}_j = \{n + j; n \in \mathcal{M}_0\}$. Therefore (1.1) follows from the fact that

$$p^a \lim_{N \to \infty} \frac{1}{N} \sum_{n=0}^{N-1} \mu\left((\phi^n \mathbf{x})_0^{m-1} = \mathbf{g}\right)$$
$$= \sum_{0 \leq j \leq p^a} \frac{1}{|\mathcal{M} \cap \{0, \ldots, N-1\}|} \sum_{n \in \mathcal{M}_j \cap \{0, \ldots, N-1\}} \mu\left((\phi^n \mathbf{x})_0^{m-1} = \mathbf{g}\right)$$
$$= \nu(\mathbf{g}) \qquad \square$$